\documentclass{amsart}
\usepackage{latexsym}

\DeclareUnicodeCharacter{200B}{}   
\DeclareUnicodeCharacter{00A0}{~}  
\DeclareUnicodeCharacter{2009}{\,} 
\DeclareUnicodeCharacter{2212}{-}  
\DeclareUnicodeCharacter{0415}{E}

\newtheorem{problem}{Problem}[section]

\def\defemb#1#2{\expandafter\def\csname #1\endcsname
                              {\relax\ifmmode #2\else\hbox{$#2$}\fi}}

\defemb{cF}{{\cal F}}
\defemb{cG}{{\cal G}}
\defemb{cI}{{\cal I}}
\defemb{cU}{{\cal U}}

\newcommand{\mL}{{\mathfrak L}}
\newcommand{\mR}{{\mathfrak R}}
\newcommand{\mJ}{{\mathfrak J}}

\newfont{\cyr}{wncyr8}
\newfont{\cyi}{wncyi8}

\begin{document}
    \baselineskip 18pt

\title{The Lyapin's notebook:\\
a collection of unsolved problems in Semigroup Theory}  

 \author{ Edited by S. Bershadsky, S. Kublanovsky, G. Mashevitzky}

\begin{abstract}
This collection presents a selected set of unsolved problems in semigroup theory, a fundamental branch of modern algebra. The publication is dedicated to the 110th anniversary of the birth of Evgeniy Sergeevich Lyapin (1914–2005), one of the founders of the field, the creator of a large, internationally recognized scientific school and the author of the world's first monograph on semigroups. 

The problems included are rooted in Lyapin's research legacy and in the subsequent work of his students. The collection covers several major directions of contemporary research: potential properties and embeddability of semigroups; structural problems and finiteness conditions in varieties; endomorphisms; solvable and unsolvable classes of finite semigroups and groups; power semigroups; inclusive varieties; and the theory of partial groupoids (pargoids). This "Notebook" serves both as a tribute to Lyapin's memory and as a roadmap for current and future research in algebraic systems. 
\end{abstract}
  
	\maketitle

\section*{Content} 
\begin{enumerate}
	\item Potential properties of semigroups, embeddability
	\item Varieties: structure problems and finiteness conditions
	\item Endomorphisms
	\item Solvable and unsolvable classes of finite semigroups and groups
	\item Power semigroups
	\item Inclusive varieties
	\item Partial groupoids (pargoids) and related models
	\item Short biography of Evgeniy Sergeevich Lyapin
	\end{enumerate}

\section*{Introduction} 

Evgeniy Sergeevich Lyapin (1914-2005) is a world famous scientist and one of the founders of Semigroup Theory - an important branch of modern algebra. He is the author of the first monograph in this field to appear in the world literature. The central theme of E.S. Lyapin's research was the composition of transformations, one of the fundamental operations in mathematics. 

This mathematical object plays a crucial role in both the exact and the human sciences, including linguistics, as well as in information technology. Lyapin's research problems and investigations continue to attract interest in both classical areas and in new mathematical theories that have emerged in recent decades.

During the 1970s–1990s, he worked successfully in another branch of algebra — the theory of partial algebraic operations. Together with his student, Professor A. E. Evseev, he published the first monograph in the world literature devoted to this subject.

E. S. Lyapin was also an outstanding teacher and the author of influential textbooks and problem books used by several generations of students. He selected problems for his students with great care, from both scientific and methodological perspectives. More than forty of his graduate students defended their dissertations under his supervision and went on to become well-known specialists.

This notebook presents a number of open problems posed by mathematics - graduates of the Lyapin's School. We dedicate this collection to the 110th anniversary of the founder of this school, as a tribute to his memory.

\section{Potential properties of semigroups, embeddability}

E.S. Lyapin introduced the concept of potential properties of semigroups in his book \cite{lyapin1}.

We say that a subset $A\subset S\times S$ (for some semigroup $S$) is potentially connected by some predicate $p\left({x,y}\right)$ in the class of semigroups $K$ if there exists an oversemigroup $T\in K$ of $S$ such that $p({a,b})$ is true for any pair $({a,b})\in A$. In this case we say that the subset $A$ is potentially $p$ connected.

Lyapin considered in his book the well-known Green's relations $\mL$ and $\mR$ as predicates and found conditions for the potential $\mL$ and $\mR$ connectedness of subsets in the class of all semigroups.

 The concept of potential properties easily extends to predicates of any arity, in particular, to the unary predicate of  invertibility.

In his work \cite{lyapin2} E.S. Lyapin finds necessary and sufficient conditions for the potential invertibility of an arbitrary element of a semigroup. Lyapin's criterion is algorithmically verifiable for finite semigroups.

In this regard, it is worth mentioning the work of E.G. Shutov \cite{Shutov} (a student of Lyapin), who proved that if any element of a semigroup $S$ is potentially invertible, then $S$ is embeddable to a group. This result makes it possible to find a system of new quasi-identities, different from the quasi-identities found by A.I. Malcev \cite{Mal'tsev}, and to resolve a number of open questions in this area.

  The paper \cite{Hall} contains a proof of the well-known result of S.I. Kublanovsky (a student of Prof. Lesohin, himself a student of Lyapin) that the problem of embedding a finite semigroup into a finite completely 0-simple one is algorithmically unsolvable.
 This result implies a negative solution to Lyapin's problem concerning the possibility of algorithmically verifying potential properties for Green’s relation $\mJ$ in finite semigroups.

In 1966, J. Rhodes posed the following problem:
Is there an algorithm that, for an arbitrary finite semigroup $S$ and for arbitrary two elements $a,b\in S$, decides whether $S$ can be embedded into a finite semigroup $T$ in which the equations $a=xby$ and $b=x'ay'$ are solvable, i.e. $a$ and $b$ divide each other?

 Based on the Kublanovsky result cited above, a negative solution to this problem was found by S. Kublanovsky and M. Sapir in \cite{KublSapir}. This led to the following problem:

\begin{problem} (S. Kublanovsky, M. Sapir)
  Describe the types of equations or systems of equations for which there exists an algorithm that decides for any values ​​of the constants whether the given system is solvable in some finite extension of the original semigroup.
\end{problem}
  
	Potential properties of semigroups include embedding problems into various recursive classes of finite semigroups. A homomorphic image of a subsemigroup $T$ of a semigroup $S$ is called a divisor of $S$ (or factor of $S$).

\begin{problem} (S. Kublanovsky)
Describe the divisors of the following subclasses of finite 0-simple semigroups:
\begin{enumerate}
	\item the class of finite Brandt semigroups;
\item the class of finite 0-simple semigroups presented in the Rees representation by a matrix that consists only of zeros and ones.
 \end{enumerate}
\noindent Divisors of the class of all finite 0-simple semigroups are described by positive universal formulas (see \cite{Kubl}).
\end{problem}

\begin{problem} (S. Bershadsky, S. Kublanovsky)
Is there an algorithm that decides whether a given finite semigroup is embeddable into a power semigroup over some finite group (periodic group) (see Chapter 5 of this notebook)?

Not every finite semigroup has this property, but every semigroup is a divisor of the power semigroup over a periodic group (see \cite{BersKubl}). For divisors of power semigroups over finite groups the answer is positive (see \cite{Henckell}).
\end{problem}

\section{Varieties: structure problems and finiteness conditions}

A variety is a class of algebraic systems that can be defined by a set (finite or infinite) of identities. Theory of varieties originates from the classical work of G. Birkhoff \cite{Birkh}, published in 1935. E.S. Lyapin and his students paid special attention to problems concerning varieties of semigroups, and the results of their pioneering work had a significant influence on subsequent research in this area.

The first examples of non-finitely based semigroup varieties were obtained by A.P. Biryukov (a student of Lyapin) in 1965 (see \cite{Biruk1}). In 1970 in his work \cite{Biruk2} he completely described the lattice of varieties of idempotent semigroups, which plays a special role in the theory of semigroup varieties. This result was independently obtained by J.A. Gerhard in \cite{G} and C.F Fennemore in \cite{F}. 

In 1969, in the paper \cite{lyapin}, E.S. Lyapin discovered a new method for constructing infinite irreducible sets of balanced identities.

The first result giving a complete description of all covers for important varieties of semigroups, together with their basis and the number of covers was obtained by A.Ya.Aizenshtat (a student of Lyapin) in 1972 in \cite{Aisenstat}, see also \cite{Aisen2}.

The first example of a finite simple semigroup that does not have a finite basis of identities, was constructed by G.I. Mashevitsky (a student of Lyapin) in \cite{Ma1}. He also constructed the first example of a finite 0-simple semigroup that does not have an irreducible basis of identities in the class of completely 0-simple semigroups (see \cite{Ma2}). It is easy to see that this example does not have an irreducible basis of identities in the class of all semigroups.

Varieties of semigroups that are finitely approximated relative to certain basic predicates were described by 
S.I. Kublanovsky (student of Prof. M.M. Lesohin, himself a student of Lyapin) in 1980 in \cite{Kubl1}. This result provided a new series of non-locally finite varieties, in which the problems associated with the corresponding predicates - such as predicates of divisibility, entry into ideals, entry into maximal subgroup - are algorithmically solvable.

It should be noted that a complete description of varieties consisting of residually finite semigroups, obtained in the cited work, was independently derived in the works of E.A. Golubov and M.V. Sapir \cite{Golubov} and R. Mackenzie \cite{McKen1}, \cite{McKen2}. Recall that an algebraic system $A$ is said to be residually finite if for any two distinct elements $a, b\in A$ there exists a finite system $F$ and a homomorphism $f: A\to F$ such that $f(a) \neq f(b)$.

It should also be noted that for many well-known predicates, S.I. Kublanovsky obtained an algorithmically verifiable descriptions in the language of indicator systems (see \cite{Kubl2}, \cite{Shevrin}).

The first examples of semigroup varieties with the undecidable equality and divisibility problems were obtained by I.L. Melnichuk (a student of Lyapin) in \cite{Melnich}.

The next natural problem was posed by E.S. Lyapin in the early 1970s. Its solution may have useful applications, for example, in automata theory.

\begin{problem} (E.S. Lyapin)
Describe the identities of the full transformation semigroup $T(n)$ on an $n$-element set. \\
A description of some classes of identities was obtained in the works \cite{PSSSV}, \cite{Ma3}.
\end{problem}

We say that a variety satisfies a property \textit{locally} if any finitely generated semigroup in the variety satisfies this property. The following problem naturally arises.

\begin{problem} (S. Kublanovsky)
Describe locally residually finite varieties of semigroups for known predicates.\\
In the global case, such description is obtained in \cite{Kubl1}, \cite{Kubl2}. For the equality predicate in the local case, the description is obtained in \cite{Sapir}. Locally residually finite varieties of rings were described in \cite{Kubl3}.
\end{problem}

Let $RS_n$ denote the variety generated by all $0$-simple semigroups satisfying the identity 
$x^2 = x^{n+2}$. Such varieties and their subvarieties are called \textit{Rees-Sushkevich varieties} (see \cite{Kubl4}).

\begin{problem} (S. Kublanovsky)
Describe the covers of the semigroup variety $RS_n$ in the lattice of varieties of all semigroups and in the lattice of varieties satisfying the identity $x^2 = x^{n+2}$. \\
It follows from the results of \cite{Kubl5} that in the latter lattice the number of covers is finite. It is also known that every subvariety of $RS_1$ has only finitely many covers (see \cite{Kubl4}).
\end{problem}

\begin{problem} (S. Kublanovsky)
Describe periodic varieties generated by regular semigroups. \\
For periodic varieties generated by completely 0-simple semigroups, such a description was obtained in \cite{Kubl4}.
\end{problem}

A word $u$ is called \textit{blocking} (respectively, \textit{non-blocking} or \textit{avoidable}) if, for every finite alphabet the set of words over this alphabet that do not contain images of $u$ under substitutions, is finite (respectively, infinite). These concepts play an important role in the theory of semigroup varieties and in theoretical combinatorics (see \cite{Bean}, \cite{Zimin}).

\begin{problem} (I. Melnichuk)
Let $c(u)$ denote the number of distinct letters occuring in the word $u$. Let $m(u)$ denote the minimal number of letters in a finite alphabet $A$ for which there exists an infinite word over $A$ that does not contain images of $u$ under substitutions. Determine the number $m(u)$ for a non-blocking word $u$ containing $c(u)$ letters. \\
Some estimates were obtained in works [\cite{Melnich1}, \cite{Melnich2}, \cite{McNulty}, \cite{Baker}].
\end{problem}

A variation of Tarski's problem on the recursiveness of the set of finite semigroups with a finite basis of identities (the finitely based semigroups) is probably the most famous and intriguing open problem in the theory of semigroup varieties. The following problem represents a special case of this question.

\begin{problem} (G. Mashevitzky)
Is the set of finite finitely based completely 0-simple semigroups over Abelian groups, in particular, over a two-element group, recursive?\\
An infinite chain of varieties, generated by completely 0-simple semigroups over the two-element group, in which finitely based and non-finitely based varieties alternate, was constructed in \cite{M4}.
\end{problem}

In the following series of problems we consider varieties of ordered semigroups, that is, semigroups equipped with a partial order compatible with the semigroup operation. Such varieties are defined by identities together with identical inequalities. Identical inequalities are often called \textit{identical relations of precedence}.
For example, the semigroup of real numbers on the segment $S = [0,1]$ with multiplication satisfies the identical inequality $xy \leq x$.

E.S. Lyapin proposed the problem of studying varieties defined by identities and identical inequalities in the late 1970s.

\begin{problem} (E.S. Lyapin)
Under what conditions is a system of identical inequalities (identical relations of precedence) equivalent to a system of identities?\\
This problem is solved for systems consisting of a single identical inequality in \cite{Mak1}.
\end{problem}

\begin{problem}  (E.S. Lyapin)
Describe all varieties of ordered semigroups that have covers in the lattice of varieties of ordered semigroups. \\
Some such varieties are described in works \cite{Mak2}, \cite{Mak3}. In the lattice of varieties of ordered semigroups, not every variety has a cover (see \cite{Mak2}).
\end{problem}

\begin{problem} (V. Makaridina)
Describe all varieties of ordered semigroups in which every semigroup has only the trivial order.\\
Some examples are given in \cite{Mak4}.
\end{problem}

\section{Endomorphisms}

An important characteristic of mathematical structures is their endomorphisms, that is transformations that preserve the structure. The set of all endomorphisms forms a semigroup under the operation of composition. Lyapin and his students paid special attention to the study of this semigroup. 

In the work \cite{lyapin3}, published in 1966, Lyapin obtained an abstract characterization of the semigroup of endomorphisms for a fairly general class of mathematical structures.

In 1969, in \cite{DlabNeum}, V. Dlab and B. Neumann constructed arbitrarily large finite semigroups with only eight endomorphisms. Their question of whether a smaller number could be achieved was answered by S. Kublanovsky (MR734514, 1983), who showed that for semigroups containing at least four elements the minimal possible number of endomorphisms is four. Moreover, there exist infinitely many finite semigroups with exactly four endomorphisms.

Let $K$ be a class of semigroups. For each natural number $n$, denote by $f_K(n)$ the minimal number of endomorphisms among all finite semigroups from $K$ having at least $n$ elements. The function $f_K(n)$ is called the \textit{endomorphic function} of the class $K$. A class $K$ is called \textit{endomorphic} if its endomorphic function is unbounded. In particular, we are interested in endomorphic varieties of semigroups. The following problems arise concerning endomorphic functions.

\begin{problem} (S. Kublanovsky)
Describe endomorphic varieties of semigroups.\\
Some sufficient conditions were obtained in \cite{kubl}.
\end{problem}

\begin{problem} (S. Kublanovsky)
Is it true that for every variety of commutative semigroups its endomorphic function is identical? \\
The variety of commutative semigroups is endomorphic (see \cite{kubl}).
\end{problem}

\begin{problem} (S. Kublanovsky)
Is it true that for any endomorphic variety of semigroups its endomorphic function is identical?
\end{problem}

\begin{problem} (S. Kublanovsky)
Obtain estimates for the values ​​of the endomorphic function for various basic classes of semigroups.
\end{problem}

In 1969, M.M. Lesokhin (a student of E. S. Lyapin) posed the problem of describing commutative semigroups that have no endomorphisms other than power functions, i.e. such endomorphisms that map each element to a fixed power of this element (\cite{Sverd}, problem 1.34). Below we call commutative semigroups whose endomorphisms are all power functions \textit{Lesokhin semigroups}.

The first nontrivial example of a Lesohin semigroup was constructed by B.K Boguta (a student of Prof. Lesokhin) in the paper \cite{Boguta}. Lesokhin's problem in the class of regular semigroups was solved by A.I. Kuptsov (a student of Prof. Lesokhin) in \cite{Kupcov}.

It can be verified by direct calculations that a finitely generated Lesokhin semigroup is either finite with a single idempotent or has no idempotents. Finite Lesokhin semigroups and $2$-generated Lesokhin semigroups without idempotents were described by A.A. Borisov in \cite{Borisov1} and \cite{Borisov2}.
The next series of problems was set by A.A. Borisov (a student of Lyapin).

\begin{problem} (A.A. Borisov)
Describe finitely generated  Lesokhin semigroups without idempotents.\\  
It is proved in \cite{Borisov4} that finitely generated separative Lesokhin semigroups without idempotents should be isomorphic to subsemigroups of an infinite monogenic semigroup.
\end{problem}

\begin{problem} (A.A. Borisov)
Describe commutative semigroups all of whose endomorphisms are injective. \\
It is proved in \cite{Borisov4} that finitely generated commutative semigroups all of whose endomorphisms are injective are precisely subsemigroups of an infinite monogenic semigroups or finite ordinal sums of such semigroups.
\end{problem}

A commutative semigroup is called an \textit{E-semigroup} if it is isomorphic to the semigroup of all its endomorphisms under pointwise multiplication. For Abelian groups, this concept was introduced in \cite{DlabNeum}. In \cite{Borisov4}, A.A. Borisov described several classes of finitely generated commutative E-semigroups: finite ones; semigroups without idempotents; cancellative semigroups; and Archimedean semigroups.

\begin{problem} (A.A. Borisov)
Describe finitely generated commutative separative E-semigroups. \\
It is proved in \cite{Borisov4} that for a finitely defined commutative semigroup the problem of deciding whether 
the semigroup of all its endomorphisms (under pointwise multiplication) is Archimedean is algorithmically solvable.
\end{problem}

\begin{problem} (A.A. Borisov)
Describe commutative semigroups whose semigroups of endomorphisms (with pointwise multiplication) is cancelative.
\end{problem}

\begin{problem} (A.A. Borisov)
Describe commutative semigroups whose semigroups of endomorphisms (with pointwise multiplication) are isomorphic to one of the following semigroups: 
\begin{enumerate}
	\item the additive semigroup of all positive rational numbers;
	\item the additive semigroup of all non negative rational numbers.
\end{enumerate}
\end{problem}

\begin{problem} (A.A. Borisov)
Describe commutative semigroups whose semigroups of endomorphisms (with pointwise multiplication) are free semigroups.
\end{problem}

It is proved in \cite{Borisov4} that for a finitely generated commutative semigroup $S$ the following
conditions are equivalent:
\begin{enumerate}
	\item the semigroup of all endomorphisms of $S$ under pointwise multiplication is a free commutative semigroup; 
	\item the semigroup of all endomorphisms of $S$ under pointwise multiplication is an infinite monogenic semigroup generated by an identical endomorphism.
\end{enumerate}
In other words $S$ is a Lesokhin semigroup.

Let $S$ be a semigroup without idempotents. For elements $a, b\in S$ we write $a < b$ if $a = ab = ba$. It is easy to verify that $<$ is a strict partial order on $S$. Hence one can speak about comparable and incomparable elements and about chains in $S$.  

The length of a finite chain is the number of its elements, and each element forms a chain of length $1$. If $S$ contains a chain of length $r$ ($r\geq 1$) but no chain of length $r+1$, we say that the \textit{chain rank} of $S$ equals $r$. If $S$ contains chains of length $n$ for all natural numbers $n$, then the chain rank of $S$ is said to be infinite.

\begin{problem} (A.A. Borisov)
Is it true that the chain rank of any finitely generated (or finitely presented) semigroup without idempotents is finite?\\
It is proved in \cite{Borisov4} that finitely generated commutative semigroups without idempotents have finite chain rank.
\end{problem}

\begin{problem} (A.A. Borisov)
Describe commutative semigroups that have exactly four endomorphisms. \\
It is proved in \cite{Borisov3} that every finite commutative semigroup with more than four elements has at least five endomorphisms.
\end{problem}

In his book \emph{Semigroups}, E.S. Lyapin emphasized that establishing connections between the properties of endomorphism semigroups and the properties of the underlying structured sets is one of the most important directions in semigroup theory. This direction may be viewed as a natural generalization of the ideas of Galois theory.

\begin{problem} (E.S. Lyapin)
For natural classes of mathematical structures - such as ordered sets, relatively free algebras, and for subsemigroups of endomorphism semigroups of these structures - determine:
\begin{enumerate}
	\item  Subclasses of structures defined up to isomorphism by these semigroups;
\item Subsemigroups of endomorphism semigroups that determine these structures up to isomorphism.
\end{enumerate}
For various classes of partially ordered sets E.S. Lyapin found endomorphism semigroups that determine these classes up to isomorphism (see \cite{Lyapin4}, \cite{L5}, \cite{L6}).
\end{problem}

For a mathematical structure $A$, denote by $End(A)$ (respectively $Aut(A)$) the semigroup (respectively group) of its endomorphisms (automorphisms) under composition. The sequence 
\[
End(A), End(End(A)), \ldots
\quad
(\text{respectively } Aut(A), Aut(Aut(A)), \ldots)
\]
is called the \textit{tower of endomorphisms} (respectively \textit{tower of automorphisms}) of $A$.

\begin{problem} (G. Mashevitzky)
Describe the endomorphism semigroup and the tower of endomorphisms (automorphisms) 
\begin{enumerate}
	\item for endomorphism semigroups of semigroups free in basic classes of semigroups, in particular in the classes of idempotent, commutative Clifford, and Clifford semigroups;
\item for endomorphism semigroups of algebras free in basic classes of mathematical structures.
\end{enumerate}
\noindent For the class of all semigroups and inverse semigroups, the description of semigroups of endo-
morphisms were obtained in \cite{MS} and \cite{MSZ}. Towers of automorphisms of free inverse
semigroups are described in \cite{Mash3}.
\end{problem}

\begin{problem} (G. Mashevitzky)
\begin{enumerate}
	\item Describe transformation semigroups that are closed in the pointwise topology.
\item Describe endomorphism semigroups that are closed in a topology of Zarisky type.
\end{enumerate}
Motivation and definitions can be found in \cite{Cameron}, \cite{Mash1}, \cite{MPP}, \cite{Mash2}. Transformation semigroups closed in the pointwise topology are precisely the semigroups of endomorphisms of mathematical structures (see \cite{Mash2}). It was shown in \cite{Mash1} that the Zariski-type topology is weaker than the pointwise topology and conditions were obtained under which these topologies coincide. Several special cases of the problem were solved in \cite{Mash2}.
\end{problem}

\section{Solvable and unsolvable classes of finite semigroups and groups}

A class of finite semigroups $\Phi$ is called recursive (or decidable) if there exists an algorithm that determines whether a given finite semigroup belongs to this class.

Let $\mathbb{H}, \mathbb{S}, \mathbb{P}$ denote the operators of taking homomorphic images, subsystems and finite direct products, respectively. These operators and their combinations are called Birkhoff operators (see \cite{Kubl6}). If $ \mathbb{K}$ is a Birkhoff operator, then $\mathbb{K}(\Phi)=\{\mathbb{K}(T)|T\in\Phi\}$. For example:\\
$\mathbb{S}(\Phi)$ is the class of semigroups embeddable into semigroups from $\Phi$;\\
$\mathbb{H}(\Phi)$ is the class of semigroups that are homomorphic images of semigroups from $\Phi$;\\
$\mathbb{H}\mathbb{S}(\Phi)$ is the class of divisors of semigroups from $\Phi$;\\
$\mathbb{H}\mathbb{S}\mathbb{P}(\Phi)$ is the pseudovariety generated by $\Phi$.

In general, the recursiveness of a class of semigroups or groups $\Phi$ does not imply the recursiveness of the class $ \mathbb{K}(\Phi)$ (see \cite{Kubl6}, \cite{Albert}, \cite{Hall1}). This leads to the following problem.

\begin{problem} (S.I. Kublanovsky)
Let $C_0$ be a class of finite $0$-simple semigroups. 
Is the class $ \mathbb{K}(C_0)$ decidable, where $\mathbb{K}$ is one of the operators $\mathbb{H}\mathbb{P}, \mathbb{H} \mathbb{P}\mathbb{S}, \mathbb{S}\mathbb{H}\mathbb{P}\mathbb{S}$?\\
For all other Birkhoff operators, the answer has been obtained in \cite{Kubl6}).
\end{problem}

There exist recursive classes $\Phi$ for which $\mathbb{K}(\Phi)$ is recursive for every Birkhoff operator $\mathbb{K}$. Such classes are called completely solvable.

From the results of \cite{Ash} and \cite{Shein} it follows that the class of finite inverse semigroups is completely solvable. Other examples of completely solvable classes are presented in \cite{Kubl6}.
In this regard, the following problems arise.

\begin{problem} (S.I. Kublanovsky)
\begin{enumerate}
	\item Is the class of finite $0$-simple semigroups over groups, belonging to a pseudovariety of groups with
decidable universal theory, completely solvable?
\item The same question for the class of finite Brandt semigroups, under a similar restriction on maximal
subgroups.
\end{enumerate}
\end{problem}

\begin{problem} (S.I. Kublanovsky)
Is the class of finite regular aperiodic semigroups completely solvable?

Recall that a semigroup is aperiodic if all its subgroups are trivial.
\end{problem}

\begin{problem} (S.I. Kublanovsky)
\begin{enumerate}
	\item Is the class of finite regular orthodox semigroups completely solvable?
	\item The same question for the class of finite regular semigroups, whose idempotents belong to a given variety of idempotent semigroups.
\end{enumerate}
\end{problem}

\begin{problem} (S.I. Kublanovsky)
Is the class of finite regular aperiodic semigroups completely solvable?

Recall that a semigroup is aperiodic if all its subgroups are trivial.
\end{problem}

The subset $\{M\in \mathbb{K}(\Phi)|\mathbb{H}\mathbb{S}\mathbb{P}(M)\subseteq\mathbb{K}(\Phi)\}$ of the class $\mathbb{K}(\Phi)$ is called the kernel of the class $\mathbb{K}(\Phi)$ and is denoted by $Cor(\mathbb{K}(\Phi))$.

\begin{problem} (S.I. Kublanovsky)
Let $\mathbb{K}$ be one of the operators $\mathbb{H}\mathbb{P}, \mathbb{H}\mathbb{P}\mathbb{S}, \mathbb{S}\mathbb{H}\mathbb{P}, \mathbb{S}\mathbb{H}\mathbb{P}\mathbb{S}$, and let $\Phi$ be an arbitrary finite class of finite semigroups. Is the class $Cor(\mathbb{K}(\Phi))$ decidable? 

For all other Birkhoff operators the answer is positive (see  \cite{Kubl6}).
\end{problem}

\section{Power semigroups}

The set of nonempty subsets of a semigroup $S$, that is the power set $P(S)$, forms a semigroup under the natural  multiplication of subsets. This semigroup is called the power semigroup of $S$ and is denoted by $P(S)$.

In his monograph “Semigroups”, E.S. Lyapin emphasized the importance of studying power semigroups from various perspectives. In 1966, at the International Congress of Mathematicians in Moscow, B.M. Schein
formulated a series of open problems in semigroup theory. In particular, he posed the following question: Does  the isomorphism of power semigroups $P(S_1) \cong P(S_2)$ imply the isomorphism of the original semigroups $S_1 \cong S_2$. The same problem was independently posed by Japanese mathematician T. Tamura in \cite{Tamura1}.
It is easy to see that in the class of groups the Tamura-Schein problem has the positive solution.
Tamura conjectured that this problem has the positive solution also in the general case. However, a negative solution in the class of all semigroups was obtained by E. Mogilyanskaya (a student of Lyapin) in 1972. Namely, in \cite{Mogil} she constructed an infinite family of infinite semigroups of arbitrary cardinality that are pairwise non-isomorphic but have isomorphic power semigroups. This result is mentioned in the monograph by J. Almeida \cite{Almeida}, pp. 357–358.

\begin{problem} (E. M. Mogiljanskaja)
Does the isomorphism of power semigroups implies the isomorphism of the original semigroups in the class of finite semigroups?
\end{problem}

A systematic study of the properties of power semigroups began in 1967 with the paper \cite{Tamura2}
by T. Tamura and J. Shafer. 
Among the results in this direction we mention an unexpected theorem of by S.G. Bershadsky, (a student of prof. M. M. Lesohin, himself a student of Lyapin) proved in \cite{Bersh}: ​​“Every semigroup can be embedded into a power semigroup over a free group."
It is also known that every commutative semigroup can be embedded into the power semigroup over an Abelian group (see \cite{Trnkova}), while not every semigroup can be embedded into the power semigroup over a periodic group (see \cite{Bersh+Kubl}).

These results lead to the following problems.

\begin{problem} (S.G. Bershadsky, S.I. Kublanovsky)
Find necessary and sufficient conditions for the embeddability of a regular semigroup into the power one over a periodic group. 

For completely regular semigroups such conditions are inversion and periodicity (see \cite{Bersh+Kubl}).
\end{problem}

\begin{problem} (S.G. Bershadsky, S.I. Kublanovsky)
Is every periodic commutative semigroup embeddable into the power semigroup over a periodic Abelian group? 

This is true for all regular commutative semigroups and for some finite non-regular commutative semigroups (see \cite{Bersh+Kubl}).
\end{problem}

\begin{problem} (S.G. Bershadsky, S.I. Kublanovsky)
Is it true that if every finitely generated subsemigroup of a semigroup $S$ is embeddable into the power semigroup over some periodic group then $S$ is also embeddable into the power semigroup over a periodic group?
\end{problem}

\begin{problem} (S.G. Bershadsky, S.I. Kublanovsky)
Find general embeddability criteria (necessary and/or sufficient conditions) for a semigroup to be embeddable into the power semigroup over a periodic group.

Some partial conditions were obtained in \cite{Donald}, \cite{Bersh+Kubl}
\end{problem}

\section{Inclusive varieties}

An inclusive variety is a class of semigroups that can be defined by universal formulas of the form $u \;\hat{\in}\; V$, where $u\in X^+$ is a word and $V\subseteq X^+$ is a finite or a countably infinite set of words. Such formulas are called \textit{identical inclusions}. If $V$ is a finite set, then the inclusive identity $u \;\hat{\in}\; V$ is a first-order formula, that is an elementary formula. An inclusive variety that cannot be defined by the first-order formulas is called a \textit{nonelementary (or non-axiomatizable) inclusive variety}.

Universal formulas of the form $U = V$, where $U, V$ are some sets of words (finite or infinite), were called \textit{collective identities} by E.S. Lyapin. He proved that every system of identical inclusions is equivalent to a system of collective identities and, conversely, every system of collective identities is equivalent to a system of identical inclusions.

These notions were introduced by E.S. Lyapin in a series of papers \cite{L1, L3, L4} in the 1970s and 1980s. It turned out that in the language of identical inclusions allows a much finer classification of various classes of groups and semigroups than the classical language of identities.

Many natural classes of semigroups, for example, nilsemigroups, nilpotent semigroups, periodic Clifford semigroups, finite groups, and cyclic $p$-groups form inclusive varieties; moreover they are nonelementary inclusive varieties. S. Bratchikov (a student of Lyapin) proved in \cite{B} that there exist continuum many nonelementary inclusive varieties of Abelian groups.

Identical inclusions constitute a special class of disjunctive identities, or infinite universally closed disjunctions. Therefore, they may find applications in the theory of varieties of group representations (\cite{PV}) and in logical geometry (\cite{Plotkin}).

Within this framework several natural problems arise. Many of them were proposed by G. Mashevitsky (a student of Lyapin).

\begin{problem} (G. Mashevitzky)
Find a theorem for inclusive varieties analogous to the Birkhoff theorem for varieties; that is, determine a (minimal) set of Birkhoff type operators such that a class of semigroups forms an inclusive variety if and only if it is closed under these operators.

Inclusive varieties are known to be closed under subsemigroups and homomorphic images.
\end{problem}

\begin{problem} (G. Mashevitzky)
Determine all inclusive varieties with finite or countable lattices of inclusive subvarieties.

A partial solution can be found in \cite{M2}.
\end{problem}

\begin{problem} (G. Mashevitzky)
\begin{enumerate}
\item Find minimal (maximal) inclusive varieties of semigroups that cannot be defined by elementary identical inclusions.
\item Describe elementary inclusive varieties of semigroups all of whose inclusive subvarieties are elementary.
\item Describe nonelementary inclusive varieties of periodic semigroups whose inclusive supervarieties in the class of periodic semigroups are nonelementary.
\end{enumerate}
All minimal inclusive varieties of Abelian groups that cannot be defined by elementary identical inclusions are described in \cite{M3}.
\end{problem}

We call two semigroups inclusively equivalent if they satisfy the same identical inclusions; equivalently if they generate the same inclusive variety.

\begin{problem} (G. Mashevitzky)
Classify semigroups from the basic classes of semigroups up to inclusive equivalence.
\end{problem}

A semigroup is said to be determined by identical inclusions up to isomorphism if every semigroup satisfying the same identical inclusions is isomorphic to it.

\begin{problem} (G. Mashevitzky)
Find semigroups from the basic classes of semigroups determined by their identical inclusions up to isomorphism. 

Two-element semilattice (\cite{L1}), quasicyclic groups (\cite{BM}) and any finite group (\cite{M2}) are determined by identical inclusions up to isomorphism.
\end{problem}

\begin{problem} (G. Mashevitzky)
Describe ordinary semigroup varieties such that their finitely generated free semigroups can be separated from each other by identical inclusions. 

The variety of semilattices satisfies this condition (\cite{M1, M01}).
\end{problem}

\begin{problem} (G. Mashevitzky)
Does there exist a nonelementary inclusive variety of idempotent semigroups?
\end{problem}

\begin{problem} (G. Mashevitzky)
Does there exist a variety of semigroups such that the class of its finite (or finitely generated) semigroups forms an elementary inclusive variety?
\end{problem}

\begin{problem} (G. Mashevitzky)
Does there exist a finite semigroup that generates an inclusive variety with an infinite lattice of inclusive subvarieties?

The lattice of inclusive subvarieties of an inclusive variety, generated by a finite group or a finite nilsemigroup is finite (see \cite{M3}).
\end{problem}


\begin{problem} (G. Mashevitzky)
Find an algorithm that decides whether the inclusive variety generated by a given finite set of finite semigroups, can be defined by a finite system of (elementary) identical inclusions or prove that such an algorithm does not exist. 
(This problem is similar to Tarski's problem for identities).

In \cite{M4}, a finite semigroup was constructed that does not have a finite basis for its elementary identical inclusions.
\end{problem}

\begin{problem} (G. Mashevitzky)
\begin{enumerate}
	\item Describe classes of semigroups in which elementary inclusive varieties form a sublattice of the lattice of inclusive varieties.
\item Describe nonelementary inclusive varieties that form an upper subsemilattice in the lattice of inclusive varieties within the class of periodic semigroups.
\end{enumerate}
\noindent Elementary inclusive varieties form a sublattice, and nonelementary inclusive varieties form an upper subsemilattice in the lattice of inclusive varieties of finite Abelian groups (\cite{M3}).
\end{problem}

\section{Partial groupoids (pargoids) and related models}

In the 1960s–1980s, E.S. Lyapin developed the theory of partial algebraic operations. The main outcome of this work was the monograph \cite{ESL1}, written jointly with his student, Prof. A.E. Evseev. 

This research was motivated by the fact that many important classes of transformations arising in different branches of mathematics form partial groupoids (called pargoids in Lyapin’s terminology) with respect to the operation of composition. For example, the class of morphisms of a category, the class of projections, and the class of closure operators form pargoids that are not semigroups.

The important notion of locally finite embeddability (that is, embeddability into a finite semigroup or group), which plays a significant role in the theory of quasivarieties and in solving word equality problems (see \cite{Evans}), is naturally formulated in the language of partial groupoids. The unsolvability of a number of well-known algorithmic problems for finite groups and semigroups has been proved using partial groupoids (see \cite{HKMST}, \cite{Kublan1}, \cite{Kublan2}).

The non-axiomatizability of certain classes of semigroups (or algebras of other signatures), for example inclusive varieties, can also be established using partial groupoids together with the Łoś–Tarski theorem (see \cite{Malcev}, \cite{GM2}).

A semigroup $S$ is called an outer (respectively, internal) semigroup extension of a pargoid $P$ if $P\subseteq S$ (respectively, $P=S$) and the restriction of the operation of $S$ to $P$ coincides with the operation in $P$.

A directed endomorphism, that is a projection of an ordered set is called a closure operator on this set.

\begin{problem} (E.S. Lyapin)
\begin{enumerate}
	\item Characterize pargoids of projections (idempotents in the semigroup of all transformations).
\item Characterize the pargoids of matrix units.
\item Characterize the pargoids of idempotents of the basic classes of semigroups.
\end{enumerate}

 Some classes of pargoids were characterized in the paper \cite{ESL5}.
\end{problem}

\begin{problem} (E.S. Lyapin)
\begin{enumerate}
	\item For the main classes of mathematical structures and natural pargoids of transformations, determine the classes of structures that are defined by these pargoids up to isomorphism.
\item Find pargoids of transformations that are determined up to isomorphism by their semigroup of endomorphisms (or group of automorphisms).
\end{enumerate}

 E.S. Lyapin proved that partially ordered sets are determined up to isomorphism by the pargoids of their closure operators on these sets (see \cite{ESL2}, \cite{ESL3}, \cite{ESL4}).
\end{problem}

\begin{problem} (E.S. Lyapin)
Characterize pargoids of transformations that admit an internal semigroup extension. 

Some sufficient conditions are given in \cite{ESL6}.
\end{problem}

A class of algebraic systems is called a class of models if its signature contains only relations. A partial binary operation can be represented by a ternary relation, which allows pargoids, groups, and semigroups to be viewed as models. This approach provides new classification possibilities, as demonstrated in the papers \cite{Gorn1}, \cite{Gorn2} by O.M. Gornostaev (a student of Prof. A.E. Evseev, himself a student of Lyapin).

A model $M$ with a ternary relation is said to be representable in a groupoid $A$ (in particular, in a semigroup or a group), if there exists an injective homomorphism from $M$ into the corresponding ternary model of $A$. This definition can be generalised to universal algebras of arbitrary signature.

If $C$ is a variety, then the class of models representable in $C$ is a quasivariety of models, that is, the representability condition can be expressed by quasiidentities (see \cite{Gorn2}).

An involution is a unary operation $x\to x^{-1}$ that satisfies the identity $(x^{-1})^{-1}=x$.
Semigroups equipped with an involution are called involuted semigroups.
Clifford semigroups and inverse semigroups are examples of involuted semigroups.

\begin{problem} (O.M. Gornostaev, S.I. Kublanovsky)
Describe varieties of involuted semigroups for which the quasivarieties of models representable in these varieties are finitely axiomatizable.

For classical varieties such a description has been obtained in \cite{Gorn1}.
\end{problem}

Short biography of Evgeniy Sergeevich Lyapin

\noindent "Cast your bread upon the waters, for you will find it after many days." - Ecclesiastes 11:1

Evgeniy Sergeevich Lyapin was an outstanding Russian mathematician, one of the pioneers of the algebraic theory of semigroups, and the author of the world’s first monograph on this subject. 
The following remarks add several details to Lyapin’s portrait that are less widely known.

A major contribution Russia made to world culture in the twentieth century was the creation of unique scientific schools. Lyapin founded one such school in algebra that became internationally recognized. Graduates of this school later trained generations of mathematics teachers at pedagogical universities across many regions of Russia and other countries.

A defining feature of Lyapin’s character was his integrity and principled conduct. Honesty, kindness, and moral steadfastness often require particular courage, and Lyapin never compromised his conscience for the sake of career advancement. He did not join political organizations—neither the Komsomol nor the Communist Party—and refused to sign collective denunciations. In difficult times he consistently defended students and colleagues who faced unjustified persecution (see some examples in Anna Ya. Aızenshtat and Boris M. Schein \textit{Evgeniy Sergeyevich Lyapin—In Memoriam}, Aequationes Math. 73 (2007) 1–9).

During World War II, Lyapin refused evacuation from besieged Leningrad together with his institute, despite repeated orders. He remained in the city and endured hunger, blockade, and bombardment. During the war he worked at the Leningrad Branch of the Main Geophysical Observatory, where he took part in calculating the load-bearing capacity of the ice road across Lake Ladoga—the so-called “Road of Life,” which saved the lives of tens of thousands of people. His work in Atmospheric Physics was considered so significant that he was awarded the Order of Lenin and granted the academic title of Professor of Atmospheric Physics.

Lyapin’s development as a scientist was inextricably linked to the traditions of the St. Petersburg mathematical school. His supervisor at the university was Professor Vladimir Abramovich Tartakovskiy, himself a student of the renowned mathematician Boris Nikolaevich Delaunay.

Lyapin defended his PhD thesis in group theory, the theory of reversible transformations. Later he turned to semigroup theory, which studies non-reversible transformations. He explained this choice simply: the vast majority of transformations in nature are not reversible.

Today the importance of semigroup theory is widely recognized. Its methods and results are used not only in classical branches of mathematics but also in automata theory, formal language theory, theoretical computer science, particle physics, and in the study of multidimensional diffusion processes in biological cells. The range of applications of semigroup theory continues to expand.

However, seventy years ago semigroup theory in the USSR was often regarded as a useless and impractical field, and serious pressure was exerted on Lyapin to abandon it. The late 1940s were a difficult period for Soviet fundamental science. By 1949 the campaign for the ideological “purity” of science had reached mathematics. At that time the university administration launched attacks on several areas, including topology, semigroup theory, and functional analysis.

Lyapin was the only participant who spoke publicly in defense of semigroup theory, mathematics and science in general, despite warnings from the administration that the decision had already been made and that he should remain silent. This episode, together with many other events from Lyapin’s life, is described in the book \textit{Mathematics across the Iron Curtain: A History of the Algebraic Theory of Semigroups} by Christopher Hollings. The attack on Lyapin and two of his colleagues at an open meeting of the Academic Council of Leningrad State University is also described in an article by N. M. Khait, \textit{The Pride of Our City}, History of Petersburg, No. 1 (23), 2005. At that meeting all those accused remained silent except Lyapin. He stated that the success of science requires the advancement of new ideas. Although history proved him right, at that time he was dismissed from Leningrad University for defending his right to work on semigroup theory. He subsequently continued his work at the Herzen State Pedagogical University of Russia (at that time Leningrad State Pedagogical Institute).

Lyapin devoted a large part of his time and talent to teaching and to improving mathematical education. In particular, he initiated important reforms in the mathematical curricula of pedagogical institutes in the Soviet Union. In the 1960s he served as an expert for UNESCO in India, where he prepared extensive recommendations for improving mathematics education.

In Lyapin’s personality the talents of a mathematician, a teacher, a writer, and a philosopher were harmoniously united. He generated many original mathematical ideas and wrote influential textbooks and monographs. At the same time, he authored several works of fiction and reflected deeply on the development of human civilization in his philosophical book \textit{Dynamics of Civilizations}. In his writings he discussed the meaning of life, the role of the scientist and the responsibilities of the teacher, and the relationship between pure and applied science. These reflections reveal a deeper understanding of the remarkable personality of Evgeniy Sergeevich Lyapin.

\end{document}